\documentclass[12pt, reqno, twoside, letterpaper]{amsart}
\usepackage{amssymb}
\usepackage{amsmath}

\usepackage[text={144mm,233mm},centering]{geometry}

\usepackage{setspace}
\usepackage{amsthm}
\theoremstyle{plain}
\newtheorem{X}{X}[section]
\newtheorem{theorem}[X]{Theorem}

\usepackage{hyperref}
\hypersetup{
    pdftoolbar=true,                        
    pdfmenubar=false,                        
    pdffitwindow=false,                      
    pdfstartview={FitH},                    
    pdftitle={A note on the theorem of Maynard and Tao},    
    pdfauthor={Tristan Freiberg},    
    pdfcreator={Tristan Freiberg},
    pdfproducer={Tristan Freiberg},
    pdfnewwindow=true,                       
    colorlinks=true,
    linkcolor={black},                      
    citecolor={black},
    filecolor={black},
    urlcolor={black},         
}

\usepackage{cite}

\numberwithin{equation}{section}

\renewcommand{\le}{\ensuremath{\leqslant}}
\renewcommand{\ge}{\ensuremath{\geqslant}}

\title[A note on the theorem of Maynard and Tao]{A note on the theorem of Maynard and Tao}

\author{Tristan Freiberg}

\address{Department of Mathematics, 
         University of Missouri, Columbia MO, USA.}
\email{freibergt@missouri.edu}

\date{\today}

\begin{document}

\begin{abstract}
As a corollary to the recent extraordinary theorem of Maynard and 
Tao, we re-prove, in a stronger form, a result of Shiu concerning 
``strings'' of consecutive, congruent primes.
\end{abstract}

\maketitle

\section{Introduction}
 \label{sec:1}

One form of the prime $k$-tuple conjecture asserts that there are 
infinitely many integers $n$ for which 
$g_1n + h_1,\ldots,g_kn+h_k$ is a $k$-tuple of primes, provided 
the $k$-tuple of linear forms  
$g_1x + h_1,\ldots,g_kx + h_k$ is {\em admissible}.
By admissible, we mean that $\prod_{i = 1}^{k}(g_ix + h_i)$ has 
no fixed prime divisor, that is for all primes $p$,
\begin{align*}
 \#\{n \bmod p : 
  {\textstyle \prod_{i = 1}^{k}(g_in + h_i) } \equiv 0 \bmod p\}
   < p,
\end{align*}
so that there is no ``obvious'' reason why there can't be 
infinitely many $k$-tuples of primes of that form.
(For inadmissible $k$-tuples, for every integer $n$, at least one 
of $g_1n+h_1,\ldots,g_kn + h_k$ is a multiple of $p$ for some 
$p$.)
For the purposes of this note, the $g_i$ shall be taken to be 
positive integers, the $h_i$ integers, and 
$g_1x + h_1,\ldots,g_kx + h_k$ shall be $k$ distinct linear 
forms.

In a stunning new development, Maynard \cite{MAY} and Tao have 
come tantalizingly close to proving this conjecture --- among the 
greatest of outstanding problems in number theory.
We refer the reader to an expository article of Granville 
\cite{GRAN} for the recent history and ideas leading up to this 
breakthrough, and a discussion of the impact it will have. 
Their incredible theorem can be stated as follows (see \cite{GRAN} 
for this formulation).

\begin{theorem}
[Maynard-Tao]
 \label{T:1.1}
For any integer $m \ge 2$, there exists an integer $k = k(m)$ 
such that the following holds.
If $g_1x + h_1,\ldots,g_kx + h_k$ is an admissible $k$-tuple, then 
for infinitely many integers $n$, there exist $m$ or more primes 
among $g_1n + h_1,\ldots,g_kn + h_k$.
\end{theorem}

Without doubt, numerous interesting applications will be made of 
this wonderful result and its proof.
Many have already been proved or alluded to by Granville 
\cite{GRAN}.
The purpose of this note is to explain the following corollary to 
the theorem of Maynard and Tao. 

\begin{theorem}
 \label{T:1.2}
Let $p_1 = 2 < p_2 = 3 < \cdots $ be the sequence of all primes.
Let $q \ge 3$ and $a$ be a coprime pair of integers, and let 
$m \ge 2$ be an integer.
There exists a constant $B = B(q,a,m)$, depending only on $q$, $a$ 
and $m$, such that the following holds.
There exist infinitely many $n$ such that 
\begin{align*}
 \text{$p_{n+1} \equiv p_{n+2} \equiv \cdots \equiv p_{n+m} 
                 \equiv a \bmod q$
       \hspace{6pt} and \hspace{6pt} $p_{n+m} - p_{n+1} \le B$.}
\end{align*}
\end{theorem}

Without the extra requirement that the ``string'' of $m$ 
consecutive, congruent primes be contained in a bounded length 
interval, the theorem was proved by Shiu \cite{MR1760689}, who 
attributed to Chowla the conjecture that there should be 
infinitely many pairs of consecutive primes $p_n$, $p_{n+1}$ such 
that $p_n \equiv p_{n+1} \equiv a \bmod q$.
Note that when we speak of consecutive primes, we mean consecutive 
terms in the sequence of all primes, not just the sequence of primes 
in the arithmetic progression $a \bmod q$, and they are not 
necessarily consecutive terms of the arithmetic progression 
$a \bmod q$.
One might refer to strings of consecutive, congruent primes as 
``Shiu strings''.

\section{Proof of Theorem 1.2}
 \label{sec:2}

The proof that there are infinitely many Shiu strings in a bounded 
length interval can be deduced almost at once from the 
Maynard-Tao theorem.
One only needs to construct an admissible $k$-tuple 
$gx + h_1,\ldots,gx + h_k$ in such a way that there can never be 
any primes ``in between'' the terms $gn + h_i$. 
 
Let $q \ge 3$ and $a$ be a coprime pair of integers, and let 
$\ell_1 < \ell_2 < \cdots $ be the sequence of all primes that are 
in the arithmetic progression $a \bmod q$.
Choose an integer $m \ge 2$ and let $k = k(m)$ be an integer that 
is large enough so that, by the theorem of Maynard and Tao, any 
admissible $k$-tuple of linear forms 
$\{g_1x + h_1,\ldots,g_kx + h_k\}$ contains at least $m$ primes 
when $x = n$, for infinitely many integers $n$.
Choose an integer $t = t(q,a,k) = t(q,a,m)$ large enough so that  
$k < \ell_{t+1}$ and $\ell_{t+k} < \ell_{t+1}^2$, which can be 
done by the prime number theorem for arithmetic progressions.%
\footnote{%
If $\pi(y;q,a)$ is the number of primes 
$\ell \le y$ with $\ell \equiv a \bmod q$, then the prime number 
theorem for arithmetic progressions gives
$\pi(y;q,a) = y/(\phi(q)\log y) + O_q(y/(\log y)^2)$.
Therefore, indeed,  
$\pi(y^2;q,a) - \pi(y;q,a) > k$ for all sufficiently large $y$, 
whence $\pi(\ell_{t+1}^2;q,a) - \pi(\ell_{t+1};q,a) > k$ for all 
sufficiently large $t$.}

Let 
\begin{align*}
 g = g(q,a,k) 
   = g(q,a,m)
   =  (\ell_{t+1}\cdots \ell_{t+k})^{-1}
      {\textstyle \prod_{p \le \ell_{t+k}} p},
\end{align*}
and consider the $k$-tuple of linear forms
\begin{align*}
 \mathcal{H} = \mathcal{H}(q,a,k) 
             = \mathcal{H}(q,a,m)
             = \{gqx + \ell_{t+1},\ldots,gqx + \ell_{t+k}\}.
%                  \{gqx + \ell_{t+1},\ldots,gqx + \ell_{t+k}\}.
\end{align*}
To see that $\mathcal{H}$ is admissible, first note that for all 
primes $p$,
\begin{align}
 \label{eq:2.1}
 \#\{\ell_{t + i} \bmod p : i = 1,\ldots,k\} < p.
\end{align}
For $k < p$ this is obvious.
For $p \le k$ note that the set on the left-hand side does not 
contain $0 \bmod p$, because we chose $t$ large enough so that 
$k < \ell_{t+i}$ (which is a prime) for each $i$.
Next, if $p \mid gq$ then $p \nmid \ell_{t+1}\cdots \ell_{t+k}$ 
because of the way we constructed $g$, and so for any integer $n$,
\begin{align*}
 (gqn + \ell_{t+1})\cdots (gqn + \ell_{t+k}) 
  \equiv \ell_{t+1}\cdots \ell_{t+k}
   \not\equiv 0 \bmod p.  
\end{align*}
If $p \nmid gq$ then 
\begin{align*}
  & \#\{n \bmod p : 
   {\textstyle 
     \prod_{i=1}^{k} (gqn + \ell_{t+i}) } 
      \equiv 0 \bmod p\} \\
    & \hspace{90pt} = 
       \#\{-(gq)^{-1}\ell_{t+i} \bmod p : i = 1,\ldots,k\}  
       < p
\end{align*}
by \eqref{eq:2.1}.

Suppose $n$ and $h$ are integers with 
$\ell_{t+1} \le h \le \ell_{t+k}$ and such that $gqn + h$ is 
prime.
Then $(gq,h) = 1$, and so, by the way we constructed $g$, $h$ is 
composed only of the primes $\ell_{t+1},\ldots,\ell_{t+k}$.
Since we chose $t$ large enough so that 
$\ell_{t+k} < \ell_{t+1}^2$, we deduce that $h = \ell_{t+i}$ for 
some $i = 1,\ldots,k$.
That is, if there are any primes at all in the interval 
$[gqn + \ell_{t+1},gqn + \ell_{t+k}]$, they must be of the form 
$gqn + \ell_{t+i}$ for some $i = 1,\ldots,k$.
Thus, after removing the composite numbers from the set 
$\{gqn + \ell_{t+1},\ldots,gqn + \ell_{t+k}\}$, the primes (if 
any) that remain are necessarily consecutive primes.
(That is, consecutive terms in the sequence of all primes.)
And they are also congruent to $a \bmod q$, because so are the 
$\ell_{t + i}$.

Indeed, by the theorem of Maynard and Tao,%
\footnote{%
In Theorem \ref{T:1.1}, the $g_i$ and $h_i$ are allowed to depend on 
$k$.
}
for infinitely many integers $n$, there are at least $m$ primes among  
$gqn + \ell_{t+1},\ldots,gqn + \ell_{t+k}$, and because of the way 
$g$ was constructed, they must be consecutive primes.
Of course, they are all also congruent to $a \bmod q$, because 
$gqn + \ell_{t+i} \equiv \ell_{t+i} \equiv a \bmod q$.
We conclude the proof by setting 
$B = B(q,a,k) = B(q,a,m) = \ell_{t+k} - \ell_{t+1}$.

In fact there exist absolute positive constants $A$ and $N$ such that 
\begin{align*}
B(q,a,k) < Aq^{N}k^2\log qk.
\end{align*}
To see this, first recall that in constructing our admissible 
$k$-tuple, once we'd chosen $k$ sufficiently large in terms of $m$, 
we only required that $t$ be large enough so that $k < \ell_{t+1}$ 
and $\ell_{t+k} < \ell_{t+1}^2$. 
Now, by Linnik's theorem (see \cite[Corollary 18.8]{MR2061214} for 
this version of it) there exists an absolute positive constant $L$ 
such that for all $y \ge q^L$, 
\begin{align*}
\pi(y;q,a) \ge \frac{Cy}{\phi(q)\sqrt{q} \log y} 
            \ge \frac{Cy}{q^{3/2}\log y},
\end{align*}
where $C$ is some absolute positive constant and 
$\pi(y;q,a)$ is the number of primes $\ell \le y$ such that 
$\ell \equiv a \bmod q$.
Thus, if $M = \max\{k,\lceil L \rceil\}$ ($\lceil L \rceil$ 
the smallest integer greater than or equal to $L$) and $D$ is 
sufficiently large, then 
\begin{align*}
\pi(Dq^{\max\{L,3/2\}}k^2\log qk;q,a) \ge M^2 + k.
\end{align*}
That is, since $\pi(\ell_{M^2+k};q,a) = M^2+k$,
\begin{align*}
\ell_{M^2+k} \le Dq^{\max\{L,3/2\}}k^2\log qk. 
\end{align*}

Now, suppose that for all $t \in \{1,2,\ldots,(k-1)M + k\}$, 
$\ell_{t+k} \ge \ell_{t+1}^2$.
Then
\begin{align*}
     \ell_{(k-1)(M+1) + k+1}
  \ge \ell_{(k-1)M + k+1}^2
   \ge \ell_{(k-1)(M-1) + k+1}^4
    \ge \cdots
     \ge \ell_{k+1}^{2(M+1)}.
\end{align*}
There is at most one $\ell_i$ in each interval of length $q$, so
 $qk < \ell_{k+1}$.
We also have $\ell_{(k-1)(M+1) + k+1} \le \ell_{M^2 + k}$.
Combining all of this gives 
\begin{align*}
(qk)^{2(M+1)} < Dq^{\max\{L,3/2\}}k^2\log qk,
\end{align*}
which is absurd if $q$ or $k$ is sufficiently large.
We conclude that there must be some 
$t \in \{1,2,\ldots,(k-1)M + k\}$
such that $\ell_{t+1} < \ell_{t+k}^2$.
For such $t$ we have
$\ell_{t+k} \le \ell_{M^2+k} < Dq^{\max\{L,3/2\}}k^2\log qk$.


\begin{thebibliography}{9}
 
 \bibitem{MR2061214}
 {Iwaniec, H.~and E.~Kowalski.}
 \href{books.google.com/books?isbn=2705613668}
      {\em Analytic number theory.}
 American Mathematical Society Colloquium Publications 53.
 American Mathematical Society, Providence, RI, 2004.

\bibitem{GRAN}
{Granville, A.}
``Primes in intervals of bounded length.'' \hspace*{\fill} \\
\url{www.dms.umontreal.ca/~andrew/CEBBrochureFinal.pdf}, 
44pp., 2013.

\bibitem{MAY}
{Maynard, J.}
``Small gaps between primes.''
(Preprint.) \url{arXiv:1311.4600}, 23pp., 2013.

\bibitem{MR1760689}
{Shiu, D.~K.~L.}
``Strings of congruent primes.''
{\em J.~London Math.~Soc.~(2)}, 61(2):359--373, 2000.



\end{thebibliography}
\end{document}